\documentclass[12pt]{article}
\usepackage{amssymb,amsmath,amscd}
\newcommand{\lto}{\longrightarrow}

\begin{document}
\begin{center}
{\Large \bf Moduli Spaces of Standard Holomorphic Bundles on a
Noncommutative Complex Torus}

\vspace{12mm} Eunsang Kim${}^{a,}$\footnote{eskim@ihanyang.ac.kr},
\ and \; Hoil Kim${}^{b,}$\footnote{hikim@knu.ac.kr} \

\vspace{5mm} ${}^a${\it
Department of Applied Mathematics,\\ Hanyang University, Ansan Kyunggi-do 425-791, Korea}\\
${}^b${\it Department of Mathematics,\\
Kyungpook  National University, Taegu 702-701, Korea}\\

\vspace{12mm}

\end{center}

\begin{center}
 {\bf Abstract}
\vspace{5mm}

\parbox{125mm}{In this paper we study the moduli space of standard
holomorphic structures on a noncommutative complex two torus. It
will be shown that the moduli space is naturally identified with
the moduli space of stable bundles on an elliptic curve. We also
propose that the mirror reflection of the noncommutative complex
torus is the mirror reflection of the elliptic curve together with
a linear foliation. From this we identify the moduli space of
super cycles on the mirror reflection with the moduli space  of
standard holomorphic bundles on a noncommutative complex torus.
 } \vfill
\end{center}
\setcounter{footnote}{0}

\pagebreak

\section{Introduction}

The noncommutative tori is known to be the most accessible
examples of noncommutative geometry \cite{Co} and the physics of
open strings \cite{SeiWi, CDS}. A noncommutative torus can be
obtained in a number of different ways. It can be thought as a
strict deformation quantization \cite{Ri0} and also can be
identified with the foliation $C^*$-algebra via Morita equivalence
\cite{Co}. Here we will choose these two approach to define
noncommutative two-tori and we will discuss how they are related
with the mirror symmetry \cite{SYZ, Tyu}.

To motivate our result, let us consider D2-brane physics on a
two-torus. In Type II string theory compactified on an elliptic
curve $X_\tau$, BPS states with fixed charge are associated with
the cohomology classes of the D-brane moduli space. Suppose that
there is a wrpped 2-brane on $X_\tau$. The characteristic classes
of the Chan-Paton bundle and the RR charge vector  are related by
the Chern character and it gives brane charges associated to a
gauge field configuration. The D-brane moduli space is identified
with the moduli space of Yang-Mills connections by the
Donaldson-Uhlenbeck-Yau Theorem. By taking T-duality, the D-brane
moduli space corresponds to the D1-brane moduli space on the dual
elliptic curve of $X_\tau$ and as shown in \cite{Tyu}, the duality
is equivalent to mirror symmetry. On the other hand, open string
theory on a D2-brane can be understood from the viewpoint of
deformation quantization \cite{SeiWi}. Thus it is described by a
noncommutative two-torus. The D-brane moduli space on $X_\tau$
corresponds to the moduli space of constant curvature connections
\cite{CDS,CoR}. Now D1-brane physics is T-dual to the D2-brane
physics. The algebra of open string field on the D1-brane is
identified with the irrational rotational $C^*$-algebra
\cite{Kaj}. This $C^*$-algebra is Morita equivalent to the
$C^*$-algebra of the Kronecker foliation on the dual torus
\cite{Co}. Thus, the deformation quantization of ordinary
two-torus is related to a foliation $C^*$-algebra by T-duality.
Since T-duality on two-tori is equivalent to mirror symmetry, it
is expected to be relevant to the mirror symmetry on
noncommutative two-tori. In other words, a deformation
quantization and Kronecker foliation on a two-torus are related by
a mirror symmetry as suggested in \cite{Fu}.

In this paper, based on the D2-brane physics, we study the
mathematical aspects of the T-duality on a noncommutative complex
torus. Holomorphic structures on a noncommutative torus has been
introduced in \cite{Sch01, DiSch}. In the same sprit of \cite{CDS,
CoR},  we show that the moduli space of stable holomorphic vector
bundles on an elliptic curve $X_\tau$ is naturally identified with
the moduli spaces of standard holomorphic structures on a basic
module on a noncommutative complex torus $T^2_{\theta,\tau}$,
whose complex structure is inherited from that of $X_\tau$. Here,
a stable bundle on $X_\tau$ is deformed to a basic projective
module on $T^2_{\theta}$ by the curvature condition. This is
discussed in section 2. In section 3, motivated by \cite{Kaj}, we
suggest a linear foliation of slope $\theta^{-1}$ on the mirror
dual of $X_\tau$ is the mirror reflection of the noncommutative
complex torus $T^2_{\theta,\tau}$. On the dual torus, the
Strominger-Yau-Zaslow fibration can be understood as a linear
foliation. We find that the foliation is deformed to a Kronecker
foliation as $X_\tau$ is deformed to a noncommutative torus. The
deformation of the foliations is seen from the suspension by the
diffeomorphisms on a circle, which is a different way of defining
Kronecker foliation on a two-torus. Finally, in section 4, we show
that the moduli space of super cycles on the mirror reflection
with the moduli space  of standard holomorphic bundles on a
noncommutative complex torus. We conclude in section 5.

\section{The moduli spaces of stable bundles
on an elliptic curve}

In this section, we show that the moduli space of holomorphic
stable bundles on an elliptic curve is naturally identified with
the moduli space of standard holomorphic bundles on a
noncommutative torus, associated to a certain topological type.

Let $X_\tau=\mathbb{C}/\mathbb{Z}+\tau \mathbb{Z}$ be an elliptic
curve whose complex structure is specified by $\tau\in\mathbb{C}$,
Im $\tau\ne 0$. For $X_\tau$, the algebraic cohomology ring is
\[A(X_\tau)=H^0(X_\tau,\mathbb{Z})\oplus
H^2(X_\tau,\mathbb{Z})\cong\mathbb{Z}\oplus\mathbb{Z}.\] The Chern
character of a holomorphic vector bundle $E$ on $X_\tau$ takes the
value in $A(X_\tau)$:
\[\text{Ch}(E)=(\text{rank }E, \text{deg }E)\in
H^0(X_\tau,\mathbb{Z})\oplus H^2(X_\tau,\mathbb{Z}),\] where
$\text{deg }E=c_1(E)=\int_{X_\tau}c_1(E)$. The slope of a vector
bundle $E$ is defined by
\[\mu(E)=\frac{\text{deg
}E}{\text{rank }E}.\] A bundle $E$ is said to be stable if, for
every proper subbundle $E'$ of $E$, $0<\text{rank }E'<\text{rank
}E$, we have
\[\mu(E')<\mu(E).\]
Every stable bundles carries a projectively flat Hermitian
connection $\nabla^E$. In other words, there is a complex 2-form
$\lambda$ on $X_\tau$ such that the curvature of $\nabla^E$ is
\[R_{\nabla^E}=\lambda\cdot \text{Id}_E\] where  $\text{Id}_E$ is
the identity endomorphism of $E$. Since
\begin{align*}
c_1(E)=\frac{i}{2\pi}\text{Tr }R_{\nabla^E}
=\frac{i}{2\pi}\lambda\cdot\text{ rank }E,
\end{align*}
we have
\[\lambda=\frac{2\pi}{i} \frac{c_1(E)}{\text{rank
}E}=\frac{2\pi}{i}\mu(E).\] Thus
\begin{equation}\label{curvature}
R_{\nabla^E}=-2\pi i \mu(E) \text{Id}_E.
\end{equation}

Let us denote by $\mathcal{M}_{n,m}^s$ the moduli space of
holomorphic stable bundles of rank $n$ and degree $m$ on $X_\tau$.
In \cite{At}, it was shown that $\mathcal{M}_{n,m}^s$ is
isomorphic to $X_\tau$ when $n$ and $m$ are relatively prime. In
other words, the points of $X_\tau$ parameterize isomorphism
classes of stable bundles of rank $n$ and degree $m$.

The pair of integers $(n,m)$ also determines the topological type
of a gauge bundle on a noncommutative torus $T^2_\theta$. For
this, we first recall some notions on vector bundles on a
noncommutative torus. A noncommutative torus $T^2_\theta$ is the
deformed algebra of smooth functions on the ordinary torus with
the deformation parameter $\theta$. The algebra is generated by
two unitaries $U_1$ and $U_2$ obeying the relation
\begin{equation}\label{nct}
U_1U_2=e^{2\pi i\theta}U_2U_1.\end{equation} The above commutation
relation defines the presentation of the involutive algebra
\[A_\theta=\{\sum_{n_1,n_2\in\mathbb{Z}^2}a_{n_1,n_2}U_1^{n_1}U_2^{n_2}
\mid a_{n_1,n_2}\in\mathcal{S}(\mathbb{Z}^2)\}\] where
$\mathcal{S}(\mathbb{Z}^2)$ is the Schwarz space of sequences with
rapid decay. According to \cite{Co1}, the algebra $A_\theta$ can
be understood as the algebra of smooth functions on $T^2_\theta$
and the vector bundles on $T^2_\theta$ correspond to finitely
generated projective (left) $A_\theta$-modules.

The infinitesimal form of the dual action of the torus $T^2$ on
$A_\theta$ defines a Lie algebra homomorphism $\delta:L\rightarrow
{\rm {Der \ }}({A}_\theta)$, where $L=\mathbb{R}^2$ is an abelian
Lie algebra and ${\rm {Der \ }}({A}_\theta)$ is the Lie algebra of
derivations of $A_\theta$. For each $X\in L$,
$\delta(X):=\delta_X$ is a derivation i.e., for $u,v\in
{A}_\theta$,
\begin{align*}
\delta_X(uv)=\delta_X(u)v+u\delta_X(v).
\end{align*}
Generators of the Lie algebra $\delta_1$, $\delta_2$ act in the
following way:
\begin{align*}
\delta_i(U_i)=2\pi i U_i \ \ \text{ and } \ \ \delta_i(U_j)=0
\text{ for }i\ne j.
\end{align*}
If $\mathcal{E}$ is a projective ${A}_\theta$-module, a connection
$\nabla$ on $\mathcal{E}$ is a linear map from $\mathcal{E}$ to
$\mathcal{E}\otimes L^*$ such that for all $X\in L$,
\begin{align}
\nabla_X(\xi u)=(\nabla_X\xi)u+\xi\delta_X(u),{\rm { \ \ \
}}\xi\in {E}, u\in {A}_\theta.\notag
\end{align}
Equivalently,
\begin{align}\label{connections}
[\nabla_i,U_j]=2\pi i\delta_{ij}\cdot U_j.
\end{align}
The curvature ${F}_\nabla$ of the connection $\nabla$ is a 2-form
on $L$ with values in the algebra of endomorphisms of
$\mathcal{E}$. That is, for $X,Y\in L$,
$${F}_\nabla(X,Y):=[\nabla_X,\nabla_Y]-\nabla_{[X,Y]}.$$
Since $L$ is abelian, we simply have
${F}_\nabla(X,Y)=[\nabla_X,\nabla_Y]$.

The Chern character of a gauge bundle on $T^2_\theta$ is an
element in the Grassmann algebra $\wedge^\cdot(L^*)$, where $L^*$
is the dual vector space of the Lie algebra $L$. Since there is a
lattice $\Gamma$ in $L$, we see that there are elements of
$\wedge^\cdot \Gamma^*$ which are integral. Now the Chern
character is the map ${\rm { \ Ch
}}:K_0({A}_\theta)\rightarrow\wedge^{\rm{ev}}(L^*)$ defined by
\begin{align}
{\rm{Ch}}(\mathcal{E})=e^{i(\theta)}\nu(\mathcal{E}).\label{Ch}
\end{align}
Here $i(\theta)$ denotes the contraction with the deformation
parameter $\theta$ regarded as an element of $\wedge^2L$ and
\[\nu(\mathcal{E})=e^{-i(\theta)}\text{Ch}(\mathcal{E})
\in\wedge^{\text{even}}(\Gamma).\] The integral element
$\nu(\mathcal{E})$ is related with the Chern character on the
elliptic curve $X_\tau$.

Now let $(n,m)\in A(X_\tau)$ be the topological type of a stable
bundle $E$ on $X_\tau$. The pair of integers defines an integral
element
$\nu(\mathcal{E})$:\[\nu(\mathcal{E})=n+mdx^{12}\in\wedge^2(L^*).\]
Then the Chern character is given by
\begin{align*}
\text{Ch}(\mathcal{E})&=e^{i(\theta)}(n+mdx^{12})\\
&=(n-m\theta)+mdx^{12}.
\end{align*}
As in the classical case, we define the slope of $\mathcal{E}$ by
the number \[\mu(\mathcal{E})=\frac{m}{n-m\theta}.\] Associated to
the curvature condition (\ref{curvature}) on the stable bundle $E$
on $X_\tau$, we define a Heisenberg commutation relation by
\begin{equation}\label{2}
F_\nabla=[\nabla_1,\nabla_2]=-2\pi i
\frac{m}{n-m\theta}=\frac{2\pi}{i}\mu(\mathcal{E}).
\end{equation}
By the Stone-von Neuman theorem, the above relation has a unique
representation. As discussed in \cite{CoR}, the representation is
just $m$-copies of the Schr\"odinger representation of the
Heisenberg Lie group $\mathbb{R}^3$ on $L^2(\mathbb{R})$, where
the product on $\mathbb{R}^3$ is given by
\[(r,s,t)\cdot(r',s',t')=(r+r',s+s',t+t'+sr').\]
Then the operators $\nabla_1$ and $\nabla_2$ are the infinitesimal
form of the representation and can be written as
\begin{align}\label{con1}
(\nabla_1f)(s)&=2\pi i\frac{ms}{n-m\theta}f(s)\\
(\nabla_2f)(s)&=\frac{df}{ds}(s)\label{con2}
\end{align}
acting on the Schwartz space
$\mathcal{S}(\mathbb{R}\times\mathbb{Z}/m\mathbb{Z})\cong
\mathcal{S}(\mathbb{R})\otimes \mathbb{C}^m$.  To specify the
module $\mathcal{E}$, we need to define a module action which is
compatible with the relation (\ref{nct}) for $T^2_\theta$. Let us
first consider unitary operators $W_1$, $W_2$ on
$\mathcal{S}(\mathbb{Z}/m\mathbb{Z})=\mathbb{C}^m$ defined by
\begin{align*}
W_1f(\alpha)&=f(\alpha-n)\\
W_2f(\alpha)&=e^{-2\pi i\frac{\alpha}{m}}f(\alpha).
\end{align*}
Then
\[W_1W_2=e^{2\pi i\frac{n}{m}}W_2W_1.\]
In other words, $W_1$ and $W_2$ provide a representation of the
Heisenberg commutation relations for the finite group
$\mathbb{Z}/m\mathbb{Z}$. Associated to the representations
(\ref{con1}) and (\ref{con2}), we have Heisenberg reprersentations
$V_1$ and $V_2$ on the space $\mathcal{S}(\mathbb{R})$ as
\begin{align*}
V_1f(s)&=e^{2\pi i(\frac{n}{m}-\theta)s}f(s)\\
V_2f(s)&=f(s+1).
\end{align*}
The operators obey the relation
\[V_1V_2=e^{-2\pi i(\frac{n}{m}-\theta)}V_2V_1.\]
Finally, the operators
\begin{align}\label{3}
U_1=V_1\otimes W_1 \
\ \ \text{ and } \ \ \ U_2=V_2\otimes W_2
\end{align} acting on the
space $\mathcal{S}(\mathbb{R}\times \mathbb{Z}/m\mathbb{Z})$
satisfy the relation
\[U_1U_2=e^{2\pi i\theta}U_2U_1.\]
Thus the module $\mathcal{E}$ is the Schwartz space
$\mathcal{S}(\mathbb{R}\times \mathbb{Z}/m\mathbb{Z})$ with the
module action (\ref{3}). We write the module as
$\mathcal{E}_{n,m}(\theta)$. This module admits a constant
curvature connection whose curvature is given by (\ref{2}). There
are other connections which have the same constant curvature as
(\ref{2}). It is known that all constant curvature (\ref{2})
connection on  connection $\mathcal{E}_{n,m}(\theta)$ are given as
\begin{align*}
\nabla_1&=2\pi i(\frac{m}{n-m\theta})s+2\pi i R_1\\
\nabla_2&=\frac{d}{ds}+2\pi i R_2
\end{align*}
where $R_1,R_2\in \mathbb{R}/\mathbb{Z}$.  Hence the moduli space
of constant curvature connections on $\mathcal{E}_{n,m}(\theta)$,
specified by (\ref{2}), is the ordinary torus $T^2$. Denote by
$\mathcal{M}_{n,m}^s(\theta)$ the moduli space of constant
curvature connections on $\mathcal{E}_{n,m}(\theta)$ with the
relation (\ref{2}). Thus we see  that the moduli space
$\mathcal{M}_{n,m}^s$ is homeomorphic to the moduli space
$\mathcal{M}_{n,m}^s(\theta)$.

We now consider holomorphic structures on the moduli space
$\mathcal{M}_{n,m}^s(\theta)$. Let us fix a complex number
$\tau\in\mathbb{C}$ such that Im $\tau\ne 0$. The parameter $\tau$
defines a complex structure on $T^2_\theta$ via derivation
$\delta_\tau$ spanning $\text{Der}(A_\theta)$:
\begin{equation}\label{comp}
\delta_\tau(\sum_{(n_1,n_2)\in\mathbb{Z}^2}a_{n_1,n_2}U^{n_1}U^{n_2})=2\pi
i\sum_{(n_1,n_2)\in\mathbb{Z}^2}(n_1\tau+n_2)a_{n_1,n_2}U^{n_1}U^{n_2}.
\end{equation}
The noncommutative torus equipped with such a complex structure is
denoted by $T^2_{\theta,\tau}$. A holomorphic structure on the
module $\mathcal{E}_{n,m}(\theta)$ compatible with the complex
structure is an operator
$\overline{\nabla}:\mathcal{E}_{n,m}(\theta)\longrightarrow\mathcal{E}_{n,m}(\theta)$
such that
\[\overline{\nabla}(\xi\cdot u)=\overline{\nabla}(\xi)\cdot
u+e\cdot\delta_\tau(\xi), \ \ \ \xi\in\mathcal{E}_{n,m}(\theta),
u\in A_\theta.\] A holomorphic structure on
$\mathcal{E}_{n,m}(\theta)$ is specified by
$\overline{\nabla}=\tau\nabla_1+\nabla_2$. Along with the
connections defined (\ref{con1}, \ref{con2}), and for
$a\in\mathbb{C}$, let
\[(\overline{\nabla}_z)(f)(s,\beta)=\frac{\partial f}{\partial
s}(s,\beta)+2\pi i((\frac{m}{n-m\theta})s +z)f(s,\beta).\] Then
$\overline{\nabla}_z$ defines a standard holomorphic structure on
$\mathcal{E}_{n,m}(\theta)$. Other holomorphic structures are
determined by translations of $R_1$ and $R_2$ in (\ref{con1}) and
(\ref{con2}).  In other words, all the holomorphic structures are
determined by the complex number $z=\tau R_1+R_2$, where
$R_1,R_2\in \mathbb{R}/\mathbb{Z}$. And $z\sim z'$ if and only if
$z\equiv z'$ mod $\tau\mathbb{Z}+\mathbb{Z}$. Furthermore the
holomorphic structures are compatible with the holomorphic
structure defined by the holomorphic structure on the stable
bundle $E$ on the elliptic curve $X_\tau$. This is easily seen
when $\theta=0$. Now the moduli space of holomorphic structures on
$\mathcal{E}_{n,m}(\theta)$ is
$\mathbb{C}/\tau\mathbb{Z}+\mathbb{Z}=X_\tau$. Thus we see that
the moduli space $\mathcal{M}_{n,m}^s$ is isomorphic to the moduli
space $\mathcal{M}_{n,m}^s(\theta)$.

\section{The mirror reflection of a noncommutative complex torus}
In this section, we briefly review the topological mirror symmetry
of \cite{SYZ} for one-dimensional Calabi-Yau manifolds. We will
mainly follow the lines of \cite{Tyu}. Then we describe a
Kronecker foliation on the mirror reflection of the elliptic curve
$X_\tau$ as the mirror partner of the noncommutative complex torus
$T^2_{\theta,\tau}$.

A complex orientation of an elliptic curve $X_\tau$ is given by a
holomorphic 1-form $\Omega$, which determines a Calabi-Yau
manifold structure on $X_\tau$. A special Lagrangian cycle of
$X_\tau$ is a 1-dimensional Lagrangian submanifold $L$ such that
\[\text{Im }\Omega|_L=0 \ \ \text{ and } \ \ \text{Re
}\Omega|_L=\text{Vol}(L)\] where the volume form is determined by
the Euclidean metric on $X_\tau$. A special Lagrangian cycle is
just a closed geodesic and hence it is represented by a line with
rational slope on the universal covering space of $X_\tau$.

Let us fix a smooth decomposition
\[X_\tau=S^1_+\times S^1_-\] which induces a decomposition of
cohomology group
\[H^1(X_\tau,\mathbb{Z})\cong \mathbb{Z}_+\oplus\mathbb{Z}_-.\]
Let $[B]\in\mathbb{Z}_-$ and $[F]\in\mathbb{Z}_+$ be generators of
the cohomology group. Then the cohomology class $[F]\in
H^1(X_\tau,\mathbb{Z})$ is represented by a smooth cycle in
$X_\tau$ and is a special Lagrangian cycle. The family of special
Lagrangian cycles representing the class $[F]$ gives a smooth
fibration
\begin{equation}\label{fibration}
\pi:X_\tau\lto S_-^1:=B
\end{equation} and the base space $B=S_-^1$ is just the moduli space of
special Lagrangian cycles associated to $[F]\in
H^1(X_\tau,\mathbb{Z})$. The unitary flat connections on the
trivial line bundle $S_+^1\times\mathbb{C}\to S_+^1$ are
parameterized by $\widehat{S^1}:=\text{
Hom}(\pi_1(S^1_+),\text{U}(1))$, up to gauge equivalences. Thus we
have the dual fibration
\begin{equation}\label{dual}
\widehat{\pi}:\widehat{X}_\tau\lto B=S_-^1
\end{equation} with fibers
\[\widehat{\pi}^{-1}(b)=
\text{Hom}(\pi_1(\pi^{-1}(b),\text{U}(1))=\widehat{S^1}.\] The
dual fibration admits the section $s_0\in \widehat{X}_\tau$ with
\[s_0\cap\widehat{\pi}^{-1}(b)=1\in\text{Hom}(\pi_1(\pi^{-1}(b),
\text{U}(1)),\] so that we have a decomposition
\[\widehat{X}_\tau=\widehat{S^1}\times S_-^1.\] Hence, associated to the class
$[F]\in H^1(X_\tau,\mathbb{Z})$, the space $\widehat{X}_\tau$ is
the moduli space of special Lagrangian cycles  endowed with
unitary flat line bundles. Furthermore, $\widehat{X}_\tau$ admits
a Calabi-Yau manifold structure. In other words,
$\widehat{X}_\tau$ is the mirror reflection of $X_\tau$ in the
sense of \cite{SYZ}. Under the K\"ahler-Hodge mirror map (see
\cite{Tyu}), a complexfied K\"ahler parameter $\rho=b+ik$ defines
a complex structure on $\widehat{X}_\tau$, where $k$ is a K\"ahler
form on ${X_\tau}$ and $b$ defines a class in
$H^2(X_\tau,\mathbb{R})/H^2(X_\tau,\mathbb{Z})$. Then
$\widehat{X_\tau}$ is the elliptic curve $\mathbb{C}^*/e^{2\pi
i\rho\mathbb{Z}}$. Similarly, the modular parameter $\tau$ of
$X_\tau=\mathbb{C}^*/q^\mathbb{Z}$, $q=\exp(2\pi i\tau)$, $\text{
Im}(\tau)>0$, corresponds to a complexified K\"ahler parameter
$\widehat{\rho}$ on
 $\widehat{X}_\tau$.

Noncommutative tori can be obtained in many different ways. In
Section 2, we understood $T^2_\theta$ as a strict deformation
quantization, \cite{Ri0}, and bundles on $T^2_\theta$ were
constructed as a deformation of bundles on an ordinary torus. On
the other hand, the algebra $A_\theta$ of functions on
$T^2_\theta$ can be defined as the irrational rotation
$C^*$-algebra, \cite{Ri}. This definition will allow us to study
$T^2_\theta$ from a geometrical point of view. In fact, such a
$C^*$-algebra is obtained from a Kronecker foliation on a torus,
({\it cf.} \cite{Co}). In below, we propose that such a foliation
structure on $\widehat{X}_\tau$ define the mirror reflection of
the noncommutative complex torus $T^2_{\theta,\tau}$.

The Kronecker or linear foliation of $\widehat{X}_\tau$ associated
to the irrational number $\theta^{-1}$ is defined by the
differential equation $dy=\theta^{-1}dx$, with natural coordinates
$(x,y)$ on the flat torus determined by the symplectic form on
$\widehat{X}_\tau$. On the covering space, the leaves are
represented by straight lines with fixed slope $\theta^{-1}$ and
every closed geodesic of $\widehat{X}_\tau$ yields a compact
transversal which meets every leaf of the foliation. Such a
transversal is represented by a line with rational slope, which
means that a compact transversal is a special Lagrangian cycle in
$\widehat{X}_\tau$.

The linear foliation can be described by the {\it suspension of
diffeomorphisms}, \cite{Law}. Let
$\alpha_{\theta^{-1}}:\widehat{S^1}\lto\widehat{S^1}$ be the
diffeomorphism defined by
\[\alpha_{\theta^{-1}}(z)=\exp(2\pi i\theta^{-1})\cdot z, \ \
z\in\widehat{S^1}.\] It is the rotation through an angle $2\pi
\theta^{-1}$. Consider the product manifold $S^1\times\mathbb{R}$
with projections $p_1$ and $p_2$:
\[\begin{CD}
\widehat{S^1}\times\mathbb{R}@>p_1>>\widehat{S^1}\\
@Vp_2VV \\ \mathbb{R}
\end{CD}\]
The product manifold $\widehat{S^1}\times \mathbb{R}$ is foliated
by the leaves $p_1^{-1}(z)=\{z\}\times\mathbb{R}$ and the
foliation transverses to the fibers of
$p_2:\widehat{S^1}\times\mathbb{R}\lto \mathbb{R}$. Consider the
$\mathbb{Z}$-action on $\widehat{S^1}\times\mathbb{R}$:
\begin{equation}\label{5}
(z,b)^n:=(\alpha_{\theta^{-1}}^n(z),b+n)=(\exp(2\pi
in\theta^{-1})\cdot z, b+n), \ \ n\in\mathbb{Z}.\end{equation} The
foliation $\{\{z\}\times \mathbb{R}\}_{z\in\widehat{S^1}}$ is
preserved by the action (\ref{5}), and thus descends to a
foliation on $\widehat{X}_\tau$. Now  the quotient
$\widehat{S^1}\times_\mathbb{Z}\mathbb{R}$ carries a 1-dimensional
foliation whose leaves are the inverse images of the projection
$\bar p_1:\widehat{S^1}\times_\mathbb{Z}\mathbb{R}\lto
\widehat{S^1}/\mathbb{Z}$ under the action (\ref{5}). Note that
the quotient $\widehat{S^1}/\mathbb{Z}$ is identified with the
space of the leaves of the foliation on
$\widehat{S^1}\times_\mathbb{Z}\mathbb{R}$. On the universal
covering space of $\widehat{S^1}$, the $\mathbb{Z}$-action on
$\widehat{S^1}$ gives an identification $\theta^{-1}\sim
n\theta^{-1}$, $n\in\mathbb{Z}$ and hence
$\widehat{S^1}/\mathbb{Z}\cong
\mathbb{R}/\mathbb{Z}+\theta^{-1}\mathbb{Z}$. The leaves are
represented by straight lines of the fixed slope $\theta^{-1}$. On
the other hand, the second projection
$p_2:\widehat{S^1}\times\mathbb{R}\lto \mathbb{R}$ becomes a
fibration
\begin{equation}\label{p2}
\bar
p_2:\widehat{S^1}\times_\mathbb{Z}\mathbb{R}\lto\mathbb{R}/\mathbb{Z}\cong
S^1
\end{equation}
whose fibers are compact transversals of the linear foliation on
$\widehat{S^1}\times_\mathbb{Z}\mathbb{R}$. Now at the cost of
defining the foliation on $\widehat{X}_\tau$, the dual fibration
(\ref{dual}) is modified to the fibration (\ref{p2}) and we have
the following projections:
\[\begin{CD}
\widehat{S^1}\times_\mathbb{Z}\mathbb{R}@>\bar p_1>>\widehat{S^1}\cong
\mathbb{R}/\mathbb{Z}+\theta^{-1}\mathbb{Z}\\
@V\bar p_2VV \\ \mathbb{R}/\mathbb{Z}=S^1,
\end{CD}\]
which determines not only a foliation but also the compact
transversals for the foliation. On the other hand, other
structures on $\widehat{X}_\tau$ such as the complex structure and
the K\"ahler structure are irrelevant to define a linear foliation
on $\widehat{X}_\tau$. So, the Calabi-Yau 1-manifold
$\widehat{X}_\tau$ equipped with the $\theta^{-1}$-linear
foliation will be denoted by $\widehat{X}_{\tau,\theta^{-1}}$.
When $\theta=0$, the foliated manifold, denoted by
$\widehat{X}_{\tau,\infty}$, is just the dual fibration
(\ref{dual}) and $\widehat{X}_\tau$ is foliated by the fibers of
the fibration. Now by identifying the base space of the fibration
(\ref{p2}) with that of the fibration $\pi:X_\tau\lto S^1_-$, we
see that the foliated manifold $\widehat{X}_{\tau,\theta^{-1}}$ is
the moduli space of special Lagrangian cycles endowed with a
unitary flat line bundle and the leaves of the
$\theta^{-1}$-linear foliation. Thus, on the mirror side, the
deformation quantization $T^2_{\theta,\tau}$ of an elliptic curve
$X_\tau$ can be seen as to define a new foliation,
$\theta^{-1}$-foliation, on the foliated manifold
$\widehat{X}_{\tau,\infty}$ and the leaves of
$\widehat{X}_{\tau,\infty}$ are served as compact transversals  of
the new foliation. In other words, a linear foliation on
$\widehat{X}_{\tau}$ is uniquely determined by the deformation
parameter. In this sense, we may conclude that the foliated
manifold $\widehat{X}_{\tau,\theta^{-1}}$ is the mirror reflection
of the noncommutative complex torus $T^2_{\theta,\tau}$. Finally,
we remark from \cite{Co} that the foliation $C^*$-algebra for the
$\theta^{-1}$-linear foliation on $\widehat{X}_\tau$ is Morita
equivalent to a $C^*$-algebra defined from a compact transversal
for the foliation. In particular, if we take the line $y=0$ in
$\widehat{X}_\tau\cong \mathbb{R}^2/\mathbb{Z}\oplus\mathbb{Z}$,
then the action of leaves of the foliation on the line defines the
$\theta$  rotation $C^*$-algebra, which is known to be the
noncommutative torus $A_\theta$, \cite{Ri}. In Section 2, we
understood $T^2_\theta$ as a deformation quantization of an
ordinary torus. Thus we see that the deformation quantization and
the linear foliation $C^*$-algebra on a tours is related by the
mirror symmetry.

\section{ The moduli spaces of supercycles}
As discussed in Section 2, a stable bundle of topological type
$(n,m)$ defines a standard holomorphic bundle
$\mathcal{E}_{n,m}(\theta)$ on $T^2_{\theta,\tau}$ by deforming
the slope $\mu(E)=\frac{m}{n}$ to $\frac{m}{n-m\theta}$. In this
section, we show that the moduli space of supercycles of the slope
$\frac{m}{n}$ on $\widehat{X}_\tau$ is naturally identified with
the moduli space of standard holomorphic structures on
$\mathcal{E}_{n,m}(\theta)$.

A supercycle or a brane on $\widehat{X}_\tau$ is given by a pair
$(\mathcal{L},A)$, where $\mathcal{L}$ is a special Lagrangian
submanifold of $\widehat{X}_\tau$ and $A$ a flat connection on the
trivial line bundle $\mathcal{L}\times \mathbb{C}\to\mathcal{L}$.
A special Lagrangian cycle in
$\widehat{X}_\tau\cong\mathbb{R}^2/\mathbb{Z}\oplus\mathbb{Z}$ is
represented by a line of rational slope, so can be given by a pair
of relatively prime integers. The lines of a fixed rational slope
are parameterized by the points of interception with the line
$y=0$. Let $\mathcal{L}_{n,m}$ be a special Lagrangian submanifold
of $\widehat{X}_\tau$ given by
\begin{equation}\label{spec}
\mathcal{L}_{n,m}=\{(ns+R_1,ms)\mid s\in\mathbb{R}/\mathbb{Z}\},
\end{equation}
so that the line has slope $\frac{m}{n}$ and $x$-intercept $R_1$.
The shift of $\mathcal{L}_{n,m}$ is represented by the translation
of $R_1$. Note that a  unitary flat line bundle on
$\mathcal{L}_{n,m}$ is specified by the monodromy around the
circle. On the trivial line bundle $\mathcal{L}\times
\mathbb{C}\to\mathcal{L}$, we have a connection one-form given by
\begin{equation}\label{conne}
A=2\pi iR_2 dx, \ \ \ x\in \mathbb{R}^2, \ \
R_2\in\mathbb{R}/\mathbb{Z}
\end{equation}
so that the monodromy between points $(x_1,y_1)$ and $(x_2,y_2)$
is $\exp[2\pi i R_2(x_2-x_1)]$. Thus, the shift of connections is
represented by monodromies. Let us denote by
$\mathcal{S}\mathcal{M}_{n,m}$ the moduli space of supercycles on
$\widehat{X}_\tau$, whose slope is $\frac{m}{n}$. It was shown in
\cite{Tyu} that the moduli space $\mathcal{S}\mathcal{M}_{n,m}$ is
isomorphic to $\mathcal{M}^s_{n,m}$, the moduli space of
topological type $(n,m)$ stable bundles on $X_\tau$. In Section 2,
we have shown that the moduli space $\mathcal{M}^s_{n,m}$ is
identified with the moduli space $\mathcal{M}^s_{n,m}(\theta)$ of
standard holomorphic bundles of slope $\frac{m}{n-m\theta}$. Thus
we see that
$\mathcal{S}\mathcal{M}_{n,m}\cong\mathcal{M}^s_{n,m}(\theta)$ .
This isomorphism can also be obtained from a geometric point of
view. In other words, we construct standard holomorphic bundles on
$T^2_{\theta,\tau}$ from a supercycle on $\widehat{X}_\tau$. This
will give us a more clear picture of mirror symmetry between
$T^2_{\theta,\tau}$ and $\widehat{X}_{\tau,\theta^{-1}}$.

Under the identification
$\mathcal{M}^s_{n,m}\cong\mathcal{S}\mathcal{M}_{n,m}$, the
special Lagrangian cycle $\mathcal{L}_{n,m}$ on
$\widehat{X}_\tau$, given by (\ref{spec}), corresponds to a stable
bundle $E$ with slope $\mu(E)=\frac{m}{n}$. By the deformation of
the slope $\mu(E)$ to $\frac{m}{n-m\theta}$, we get a standard
holomorphic bundle $\mathcal{E}_{n,m}(\theta)$ on
$T^2_{\theta,\tau}$. The bundle $\mathcal{E}_{n,m}(\theta)$ can be
defined directly from the special Lagrangian cycle
$\mathcal{L}_{n,m}$. Our construction is basically based on
\cite{Co}. Let us first consider a simple case when $n=1$, $m=0$.
In this case, the special Lagrangian cycle is represented by the
line $y=0$ in $\widehat{X}_\tau\cong\mathbb{R}^2/\mathbb{Z}^2$ and
it corresponds to the trivial line bundle on $X_\tau$. On the
foliated torus $\widehat{X}_{\tau,\theta^{-1}}$, the line $y=0$ is
a compact transversal for the $\theta^{-1}$-linear foliation and
each leaf meets the line countably many points. Associated to the
intersection points with the line $y=0$, each leaf defines the
rotation through the angle $\theta$ on $S^1$, which gives a
$\mathbb{Z}$-action on $S^1$. The $C^*$-algebra of the group
action, defined from the algebra of compactly supported smooth
functions on $S^1\times \mathbb{Z}$, is known to be the irrational
rotation $C^*$-algebra or the noncommutative torus $A_\theta$ (see
\cite{Co} for details). Regarding $A_\theta$ as a free module
$\mathcal{E}_{1,0}(\theta)$ on $T^2_\theta$, connections on
$\mathcal{E}_{1,0}(\theta)$ are simply the derivations on
$A_\theta$. Together with the complex parameter $\tau$, determined
by the complexfied K\"ahler form on $\widehat{X}_\tau$, the
derivations on $A_\theta$ give a holomorphic structure on the free
module $\mathcal{E}_{1,0}(\theta)$ as given in (\ref{comp}). Thus
we see that the special Lagrangian cycle represented by the line
$y=0$ naturally defines the trivial line bundle
$\mathcal{E}_{1,0}(\theta)=A_\theta$, just as it defines the
trivial line bundle on $X_\tau$. This is one reason why we take
$\widehat{X}_{\tau,\theta^{-1}}$ instead of
$\widehat{X}_{\tau,\theta}$ as the mirror reflection of
$T^2_{\theta,\tau}$.

On the other hand, the leaves of the $\theta^{-1}$-linear
foliation rotate the special Lagrangian cycle $\mathcal{L}_{n,m}$,
$m\ne 0$, different angle $\theta'$ from $\theta$. Thus the leaf
action on $\mathcal{L}_{n,m}$ defines another noncommutative torus
$A_{\theta'}$, which is strongly Morita equivalent to $A_\theta$,
so that $A_{\theta'}\cong\text{End}_{A_\theta}(\mathcal{E})$. The
finitely generated projective $A_\theta$-module
$\mathcal{E}=\mathcal{E}_{n,m}(\theta)$ is obtained from the space
of leafwise paths starting from the line $y=0$ and ending at
$\mathcal{L}_{n,m}$. Associated to each intersection points of the
leaves and the line $\mathcal{L}_{n,m}$, the $x$-intercepts of the
leaves are parameterized by
\begin{equation}\label{equ}
x=\frac{n-m\theta}{m}t+R_1 \ \ \ \text{ modulo }1,  \ \
t\in\mathbb{R}.\end{equation} which determines the manifold
\[E_{n,m}=\{((x,0),t)\in\widehat{X_\tau}\times\mathbb{R}\mid
m(x-R_1)=(n-m\theta)t \ \ \text{ modulo }1\}.\] From the equation
(\ref{equ}), one finds that $E_{n,m}$ is the disjoint union of
$m$-copies of $\mathbb{R}$. Thus the space of compactly supported
smooth functions on $E_{n,m}$ is the Schwarz space
$\mathcal{S}(\mathbb{R}\times
\mathbb{Z}/m\mathbb{Z})\cong\mathcal{S}(\mathbb{R})\otimes
\mathbb{C}^m$. Let $W_1$ and $W_2$ be unitary operators in
$\mathbb{C}^m$ such that $W_1^m=W_2^m=1$ and $W_1W_2=\exp(2\pi
i\frac{n}{m})W_2W_1$. Also, we define operators $V_1$ and $V_2$ on
$\mathcal{S}(\mathbb{R})$, from the equation (\ref{equ}), by
\begin{align}\label{vis}
V_1f(t)&=\exp(2\pi iR_1)\exp[{2\pi i(\frac{n-m\theta}{m})t}]f(t)\\
V_2f(t)&=f(t+1), \end{align} so that $V_1V_2=\exp[-{2\pi
i(\frac{n-m\theta}{m})}]V_2V_1$. Then the action of $T^2_\theta$
on $\mathcal{E}_{n,m}(\theta)$ is determined by $U_1=V_1\otimes
W_1$ and $U_2=V_2\otimes W_2$. In other words, the space
$C_c^\infty(E_{n,m})\cong\mathcal{S}(\mathbb{R})\otimes
\mathbb{C}^m$ is the basic module $\mathcal{E}_{n,m}(\theta)$ on
$T^2_\theta$, discussed in Section 2. Now the manifold together
with the module action is uniquely determined up to translation of
$R_1$. Thus the moduli space is identified with $S^1$. To get the
full information of $\mathcal{E}_{n,m}(\theta)$, one has to take
into account of the constant curvature connections on
$\mathcal{E}_{n,m}(\theta)$. Form the definition
(\ref{connections}) of connection, constant curvature connections
are uniquely determined up to translations and we can define the
connections as in Section 2,
\begin{align*}
\nabla_1&=2\pi i(\frac{m}{n-m\theta})t+2\pi i R_1\\
\nabla_2&=\frac{d}{dt}+2\pi i R_2
\end{align*}
so that $\nabla_1$ determines the position of $E_{n,m}$ or the
special Lagrangian cycles and $\nabla_2$ determines the holonomy
or the monodromies of the connection $A$. Finally, the complexfied
K\"ahler form on $\widehat{X}_\tau$ defines the complex structure
$\nabla_1+\tau\nabla_2$ so that we see that
$\mathcal{S}\mathcal{M}_{n,m}\cong\mathcal{M}_{n,m}^s(\theta)$.

\section{Discussion}
In this paper, we proposed that a linear foliation on an ordinary
symplectic  two-torus  can be regarded as a mirror reflection of a
noncommutative complex two-torus. Analogously we have discussed
that the deformation quantization and the linear foliation
$C^*$-algebra are related by mirror symmetry. Also, we showed that
the relevant moduli spaces are naturally identified. We may extend
our study to the homological mirror symmetry of Kontsevich
\cite{Ko} based on the works \cite{PolSch, Kaj, PolZa} and a
recent work \cite{Pol}. In doing so, one may need to consider the
Floer homology for foliations which was introduced in \cite{Fu}.
In the case of 2-dimensional complex tori, it is not exactly the
noncommutative version in the sense of \cite{Fu}. However, it will
be interesting to compute explicitly the Floler homology for the
foliations discussed in this paper and compare with the tensor
product of holomorphic vectors on a noncommutative complex torus,
computed in \cite{DiSch}.

%\pagebreak
%\begin{center}
%{\large \bf Acknowledgments} \\
%\end{center}
%\vspace{1mm}

\end{document}